# On Semiprime Rings with (α,α)-Symmetric Derivations

Mehsin Jabel Atteya[*], Dalal Ibraheem Rasen

Al-Mustansiriyah University, College of Education, Department of Mathematics, Baghdad, Iraq

**Email address**

mehsinatteya@yahoo.com (M. J. Atteya), dalalresan@yahoo.com (D. I. Rasen)

**Abstract**

The main purpose of this paper is to study and investigate concerning a (α,α)-symmetric derivations D on semiprime rings and prime rings R, we give some results when R admits a (α,α)-symmetric derivations D to satisfy some conditions on R.(i)$D([x,y]^{n+1})=0$ for all x, y∈ R. (ii) $[D(x^{n+1}),α(y)] = 0$ for all x, y ∈R. (iii) $[[D(x),α(x)],α(x)]= 0$ for all x ∈R. Where α: R → R is an automorphism mapping.

## 1. Introduction

Several researchers always ask why derivation? derivations on rings help us to understand rings better and also derivations on rings can tell us about the structure of the rings. For instance a ring is commutative if and only if the only inner derivation on the ring is zero. Also derivations can be helpful for relating a ring with the set of matrices with entries in the ring (see, Pajoohesh, 2007). Derivations play a significant role in determining whether a ring is commutative, see (Herstein, 1978), (Bell, Daif, 1995) and (Andima, and Pajoohesh, 2010). Derivations can also be useful in other fields. For example, derivations play a role in the calculation of the eigenvalues of matrices (see, Baker, 1959) which is important in mathematics and other sciences, business and engineering. Derivations also are used in quantum physics (see, Da Providencia, 1994). A lot of work has been done in this field [see (Bell, and Martindale, 1987; Bresar, 1993)]. In (Lee, 2001), he was introduce the (α,α)-derivation and α-commuting mappingin the following way: If $[f(x),α(x)] = 0$ for all x ∈R, then f is said to be α-commuting, where α is an automorphism. An additive map D:R→ R is said to be an (α,α)-derivation if $D(xy)=D(x)α(y)+α(x)D(y)$. In (Wu, and Zhang, 2006), Wu and Zhan studied the properties of Jordan (α,α)-derivation. Derivations arealso generalized as α-derivations, (α,β)-derivations and have been applied in the solution of some functional equations (see, e.g., (Bresar, 1992). (Park, 2009) he was proved, let n ≥ 2 be a fixed positive integer and let R be a non-commutative n!-torsion free prime ring. Suppose that there exists a symmetric n-derivation $\Delta: R^n→R$ such that the trace δ of Δ is commuting on R. Then we have Δ = 0.(Çeven and Öztürk ,2007),them proved let N be a 2-torsion free 3-prime near-ring, D a symmetric bi-(σ,τ)-derivation of N and d the trace of D. If xd(N) = 0 for all x ∈ N, then x = 0 or D = 0,where a near ring N is 3-prime if aNb = {0} implies that a= 0 or b= 0,and a mapping D:N×N →N is said to be symmetric if D(x,y) = D(y,x) for all x,y ∈N. A mapping d:N →N denoted by d(x) = D(x,x) is called the trace of D where D:N×N→N is a symmetric mapping. It is obvious that, if D:N×N →Nis a symmetric mapping which also bi-additive (i.e., additive in both arguments), then the trace of D satisfies the relation d(x+y) = d(x)+2D(x,y)+d(y) for all x, y ∈N.A symmetric bi-additive mapping D:N×N →N is called a symmetric bi-derivation if D(xy,z) = D(x,y)y+ xD(y,z) is fulfilled for all x, y, z ∈ N. In this paper we investigate concerning a (α,α)-symmetric derivations D on semiprime rings and prime rings R.

## 2. Preliminaries

Throughout, R is an associative semiprime ring. We shall write [x,y] for xy − yx. Then [xy,z]= x[y,z] + [x,z]y; [x,yz]=y[x,z]+ [x,y]z for all x, y,z ∈ R. Recall that R is prime if aRb=(0) implies a= 0 or b = 0 and semiprime if aRa=(0) implies a = 0. An additive mapping D from R in to itself is called a derivation if D(xy) = D(x)y + xD(y) for all x, y ∈R. A mapping f of R into itself is called commuting if [f(x),x]=0,and centralizing if [f(x),x] ∈Z(R) for all x ∈ R, Z(R) denotes the centre of R. An additive mapping F from R to R is said to be a commuting(resp. centralizing) if [F(x),x]=0 (resp.[F(x),x] ∈ Z(R)) holds for all $x$ ∈ R, and d is said to be central if F(x) ∈ Z(R) holds for all $x$ ∈ R. More generally, for a positive integer n, we define a mapping F to be n-commuting if [F(x), $x^n$] = 0 for all $x$ ∈ R. An additive (α,α)-derivation D:R → R is called an (α, α)-symmetric derivation on a ring R if D(xy) = D(yx) for all x,y ∈ R. (Brešar ,and Vukman,1989) have introduced the notion of a reverse derivation as an additive mapping d from a ring R into itself satisfying d(xy) = d(y)x+ yd(x),for all x, y ∈ *R*. Obviously, if R is commutative, then both derivation and reverse derivation are the same. An additive map D:R → R is said to be an (α, α)-derivation if D(xy) = D(x)α(y)+α(x)D(y) for all *x, y* ∈ *R*.

We shall need the following well-known and frequently used lemmas.

*Lemma 2.1. (Shu-Hua, and Feng-Wen, 2002: Corollary 1)*

Let R be a semiprime ring and D derivation of R,U a nonzero left ideal of R and $r_R(U)$=0.If D is centralizing or skew-centralizing on U, then D(R) ⊆Z(R) and the ideal generated by D(R) is in the Z(R).

*Lemma 2.2*. *(Samman,and Alyamani, 2007: Proposition 3.1.)*

A mapping D on a semiprime ring R is a reverse derivation if and only if it is a central derivation.

*Lemma 2.3. (Javed et al., 2012: Proposition 3.3.)*

Let R a ring and D be the (α, α)-symmetric mapping on a ∈ R, then the following statements hold for all x, y ∈R.

(i) D([x,y]) = D([y,x])= 0
(ii) D(x[x,y])= 0 = D([x,y]x)
(iii) $D([x,y])^2$ = 0.

*Note:* If R is a ring with unity, Inv(R) is the set of all invertible elements of R, then (Inv(R),.) form a group. Let a ∈ Inv(R) then x → $a^{-1}$xa is an automorphism, referred as inner automorphism, the collection Inv(R) forms subgroup of Aut(R). We denoted to the right annihilator ideal by $r_R(U)$=0,where U a nonzero left ideal of R.

## 3. The Main Results

*Theorem 3.1.*

Let R be a 2-torsion free ring and D an (α,α)-symmetric derivation and α an automorphism such that [D(x),α(x)]=0 for all x ∈ R, then D is a reverse derivation.

*Proof:* By hypothesis [D(x),α(x)] = 0 for all x ∈R. Linearization gives

$$[D(x),\alpha(y)]=[\alpha(x),D(y)] \quad (1)$$

Then

$$[D(x),\alpha(y)]= [D(y),\alpha(x)] \quad (2)$$

Using (1) and (2), we get [D(y),α(x)]=[α(x),D(y)] implies that D(y)α(x)−α(x)D(y)= 0 this gives [D(y),α(x)] = 0 for all, x, y ∈R. Replacing x by $\alpha^{-1}$(w), we get [D(y),w] = 0 for all w, y ∈R. Hence D is a central mapping. Then according to Lemma 2.2., a mapping D on a semiprime ring R is a reverse derivation.

*Theorem 3.2.*

Let R is a ring and D is an (α,α)-symmetric mapping on a ring R, then D($[x,y]^{n+1}$) =0 for all x,y∈ R, where n is a positive integer.

*Proof:* We can prove the theorem with the help of mathematical induction**.**

(i) When n=1, then we have D($[x,y]^2$) =0**,** As D is an (α,α)-symmetric mapping on a ring R, so by D is an (α, α)-symmetric mapping on a ring R, we have D([x,y]) = D([y,x]) = 0 for all x, y ∈R. Then

D([x,y]) = 0 for all x, y ∈R. Also
D($[x,y]^2$)=D([x,y])α([x,y])+α([x,y])D([x,y]). Since, we have
D([x,y])=0 for all x, y ∈R. Therefore, D($[x,y]^2$) = 0.
(ii) When n=m, then D($[x,y]^{n+1}$) =0 for all x∈ R.
(iii) Suppose that true when n=m+1, then
D($[x,y]^{m+1}$) =0 for all x∈ R.
D$([x,y])^m$ [x,y]) =0. Then
D$([x,y])^m$α([x,y])  +α($[x,y]^m$)D([x,y])=0. Accoring to Lemma 2.3.(iii), we get

$$\alpha([x,y]^m)D([x,y])=0.$$

Again according D is an (α,α)-symmetric mapping on a ring R, we have D$([x,y])^m$ =0 for all x, y ∈R. Then by right-multiplying by α([x,y]),we get

D$([x,y])^m$ α([x,y])=0 for all x, y ∈R. Then from them relations, we get

D$([x,y])^m$ α([x,y])+ α($[x,y]^m$)D([x,y])=0 for all x, y ∈R.
D$([x,y])^m$ =0 for all x, y ∈R. Right-multiplying by D([x,y]) for all x, y ∈R, leads to D$([x,y])^{m+1}$ =0 for all x, y ∈R. This completes the proof.

*Theorem 3.3.*

Let R is a prime ring and D be a (α,α)-symmetric derivation on R, then either R is commutative or D is zero.

*Proof:* By Lemma2.3 (ii) D(x[x,y])= 0 implies that
D(x)α([x,y])  +α(x)D([x,y])= 0  also  in  view  of Lemma2.3(i), we have D(x)α([x,y]) = 0
for all x,y ∈ R. Further in view of Lemma2.3(i),D([x,y]x)=0 gives

$$\alpha([x,y])D(x) = 0.$$

Replacing y by yz in α([x,y])D(x) = 0, and using it again, we get

$$\alpha([x,y])\alpha(z)D(x) = 0 \quad (3)$$

Replacing x by x + u in α([x,y])D(x) = 0, we get

$$\alpha([x,y])D(u) = \alpha([y,u])D(x) \text{ for all } x, y, u \in R \quad (4)$$

Pre multiplying (4) by α([x, y])D(u)v we get

$$\alpha([x,y])D(u)v\alpha([x,y])D(u)= \alpha([x,y])D(u)v\alpha([y,u])D(x). \quad (5)$$

Replacing z by $\alpha^{-1}(D(u)v\alpha([y,u]))$ in (3) and then using in (5), we get

$\alpha([x,y])D(u)v\,\alpha([x,y])D(u)= 0$ for all x, y, u, v $\in$ R

As R is a semiprime, therefore, $\alpha([x,y])D(u) = 0$.

Replacing y by yz in $\alpha([x,y])D(u) = 0$, we get

$$\alpha([x,y])\alpha(z)D(u) = 0.$$

Replacing z by $\alpha^{-1}(w)$, we get $\alpha([x, y])wD(u) = 0$ for all x, y,w, u $\in$R. As R is a prime ring, therefore either $\alpha([x,y])= 0$ or $D(u) = 0$ for all x, y, u $\in$R. By other word, we have $D(u) = 0$ for some u $\in$R, or $\alpha([x,y]) = 0$. That is $([x,y]) = 0$ for all x, y $\in$R and hence R is commutative or D is zero $(\alpha,\alpha)$-symmetric derivation on R.

*Theorem 3.4.*

Let R is a semiprime ring and D be an $(\alpha,\alpha)$-symmetric derivation on R, then $[D(x^{n+1}), \alpha(y)] = 0$ for all x,y $\in$R, where n is a positive integer.

*Proof:* Since D is an $(\alpha,\alpha)$-symmetric derivation, according to Lemma 2.4, replacing z by xz in $D(x[y,z]) = 0$ for all x,y,z $\in$ R, and then using $D(x[y,z]u) = 0$ for all x,y,z $\in$ R. Again replacing z by xz in the equation obtained and using $D(x[y,z]u) = 0$ for all x,y,z $\in$ R, again, we have $D(x^3[y, z]) = 0$.

Continuing this process, we get $D(x^{n+1}[y,z]) = 0$ for all x, y, z $\in$ R.

As D is an $(\alpha,\alpha)$-symmetric derivation, therefore

$$D([y,z]x^{n+1}) = 0 \text{ for all } x, y, z \in R.$$

$D(x^{n+1}[y,z]) = 0$ for all x, y, z $\in$R. In view of Lemma 2.3 (i),
$$D(x^{n+1})\alpha([y,z]) = 0. \quad (6)$$

As D be an $(\alpha,\alpha)$-symmetric derivation, therefore, $D([y,z]x^{n+1}) = 0$. In view of Lemma 2.3 (i), we have

$$\alpha([y,z])D(x^{n+1}) = 0. \quad (7)$$

Subtracting (7) from (6), we get $[D(x^{n+1}), \alpha([y,z])] = 0$.

Replacing z by yz and used again it, we get $[D(x^{n+1}), \alpha(y)]\,\alpha([y,z]) = 0$. Aftre that

*Theorem 3.6.*

Let R be a 2-torsion free semiprime ring ,U a nonzero left ideal of R such that $r_R(U)=0$,and $\alpha$: R $\to$ R is an automorphism. If there exists an $(\alpha,\alpha)$-derivation D:R $\to$ R such that $[[D(x),\alpha(x)],\alpha(x)]= 0$ for all x$\in$U, then $D(R) \subseteq Z(R)$ and the ideal generated by D(R) is in the Z(R).

*Proof:* At first , we need to define a mapping H (.,.): R×R $\to$ R

by $H(x,y)= [D(x),\alpha(y)]+[D(y),\alpha(x)]$ for all x$\in$ U, y $\in$R.

Then it is easy to see that H(x,y) = H(y,x)
for all x$\in$ U, y, z $\in$R, and additive in both arguments.

By using the Jacobian identities and definition of H(x,y), we can conclude the following:

$$H(xy,z)=H(x,z)\alpha(y)+\alpha(x)H(y,z)+D(x)[\alpha(y),\alpha(z)]+[\alpha(x),\alpha(z)]D(y). \quad (8)$$

Define the mapping h:R$\to$ R such that h(x)= H(x,x), then

$$h(x)=2[D(x),\alpha(x)]. \quad (9)$$

One can conclude h(x+y) = h(x)+h(y)+2H(x,y).
By $[[D(x),\alpha(x)],\alpha(x)]= 0$,we have,
$[h(x),\alpha(x)]= 0$.
Linearization of last equation gives

$$[h(x),\alpha(y)]+[h(y),\alpha(x)]+2[H(x,y),\alpha(x)]+2[H(x,y),(y)]=0. \quad (10)$$

Replacing x by −x in (10) and using the fact h(−x) = h(x), we get

$$[h(x),\alpha(y)]-[h(y),\alpha(x)]+2[H(x,y),\alpha(x)]-2[H(x,y),\alpha(y)]=0. \quad (11)$$

Adding (10) and (11) and using the fact that R is 2-torsion free, we get

$$[h(x),\alpha(y)]+2[H(x,y),\alpha(x)]=0. \quad (12)$$

replacing z by $z\alpha^{-1}D(x^{n+1})$ and again replacing z by $\alpha^{-1}(w)$ $w[\alpha(y),D(x^{n+1})] = 0$. As R is semiprime, therefore, we in the result obtained, we get $[D(x^{n+1}), \alpha(y)]w[\alpha(y),D(x^{n+1})] = 0$. have $[\alpha(y),D(x^{n+1})] = 0$ implies that $[D(x^{n+1}), \alpha(y)] = 0$ for all x, y $\in$ R. This completes the proof.

*Corollary 3.5*

Let D be an $(\alpha,\alpha)$-symmetric derivation on a semiprime ring R, such that $[D(x^{n+1}),\alpha(y)] = 0$ for all x,y $\in$R, where n is a positive integer, then $D(x^{n+1})$ is n+1-$\alpha$-central (resp. n+1-$\alpha$-commuting) of R.

Replacing y by xy in (12) and using H(x,y) = H(y,x), we get

$[h(x), \alpha(x)\alpha(y)]+2[H(xy,x), \alpha(x)] = 0$

In view of (8) and (9), the last expression becomes

$\alpha(x)([h(x),\alpha(y)]+2[H(x,y),\alpha(x)])+2h(x)[\alpha(y),\alpha(x)]+2[D(x),\alpha(x)][\alpha(y),\alpha(x)]+2D(x)[[\alpha(y), \alpha(x)], \alpha(x)] =0$.

In view of (9) and (12), the last equation reduces to

$$3h(x)[\alpha(y),\alpha(x)]+2D(x)[[\alpha(y),\alpha(x)],\alpha(x)]=0. \quad (13)$$

Replacing y by yx in (12) then using (9), (12) and the fact that $[h(x), \alpha(x)] = 0$,we get

$$3[\alpha(y),\alpha(x)]h(x)+2[[\alpha(y),\alpha(x)],\alpha(x)]D(x)=0. \quad (14)$$

Replacing y by yz in (13), we get

$3h(x)[\alpha(yz), \alpha(x)]+2D(x)[[\alpha(yz),\alpha(x)]\alpha(x)]$

$= 3h(x)\alpha(y)[\alpha(z),\alpha(x)]+3h(x)[\alpha(y),\alpha(x)]\alpha(z)= 0$

In view of (13) the last equation becomes

$$3h(x)\alpha(y)[\alpha(z),\alpha(x)]+4D(x)[\alpha(y),\alpha(x)][\alpha(z),\alpha(x)]+ 2D(x)\alpha(y)[[\alpha(z),\alpha(x)],\alpha(x)]=0. \quad (15)$$

Putting y= $\alpha^{-1}(D(x))$, z = y and using (9) the last equation becomes:

$$3h(x)D(x)[\alpha(y),\alpha(x)]+2D(x)h(x)[\alpha(y),\alpha(x)]+2(D(x))2[[\alpha(y),\alpha(x)],\alpha(x)]=0 \quad (16)$$

Pre multiplying (13) by D(x) and subtracting from (16), we get

$$(3h(x)D(x)-D(x)h(x))[\alpha(y),\alpha(x)]=0. \quad (17)$$

Post multiplying (17) by $\alpha(z)$, we get

$$(3h(x)D(x)-D(x)h(x))[\alpha(y),\alpha(x)]\alpha(z)=0. \quad (18)$$

Replacing y by yz in (17), we get

$$(3h(x)D(x)-D(x)h(x))(\alpha(y)[\alpha(z),\alpha(x)]+[\alpha(y),\alpha(x)]\alpha(z))=0. \quad (19)$$

Subtracting (18) from (19), we get

$$\{3h(x)D(x)-D(x)h(x)\}\alpha(y)[\alpha(z),\alpha(x)]=0. \quad (20)$$

Replacing z by $\alpha^{-1}(2D(x))$ and y by $\alpha^{-1}(t)$ in (20), we get

$$\{3h(x)D(x)-D(x)h(x)\}th(x)=0 \text{ for all } x\in U, y,t \in R \quad (21)$$

multiplying (21) by $3D(x)$ and again replacing t by $tD(x)$ in (21) and then subtracting there results, we get

$$\{3h(x)D(x)-D(x)h(x)\}t\{3h(x)D(x)-D(x)h(x)\}=0. \quad (22)$$

Semiprimeness of R implies that $3h(x)D(x)-D(x)h(x) = 0$. That is

$$3h(x)D(x) = D(x)h(x). \quad (23)$$

Now replacing y by zy in (14) and using (14) again and then putting
$y = \alpha^{-1}(D(x))$ and $z = y$, we get
$3[\alpha(y),\alpha(x)]D(x)h(x)+2[\alpha(y),\alpha(x)]h(x)D(x)+$

$$2[[\alpha(y),\alpha(x)],\alpha(x)](D(x))2=0. \quad (24)$$

Post multiplying (14) by $D(x)$ and using in (24), we get

$$[\alpha(y),\alpha(x)](3D(x)h(x)-h(x)D(x))=0. \quad (25)$$

Replacing y by zy in (25) and using (25) again, we get
$[\alpha(z), \alpha(x)]\alpha(y)(3D(x)h(x)- h(x)D(x)) = 0$
Replacing z by $\alpha^{-1}(2D(x))$ in the last equation and using (9), we get

$$h(x)\alpha(y)(3D(x)h(x)- h(x)D(x)) = 0. \quad (26)$$

Pre multiplying (24) by $3D(x)$ and replacing y by $\alpha^{-1}(t)$, we get

$$3D(x)h(x)t(3D(x)h(x)-h(x)D(x)) = 0. \quad (27)$$

Replacing y by $\alpha^{-1}(D(x)t)$ in (26), we get

$$h(x)D(x)t(3D(x)h(x)- h(x)D(x)) = 0. \quad (28)$$

Subtracting (28) from (27) and using the fact that R is semiprime, we get

$$3D(x)h(x) = h(x)D(x). \quad (29)$$

Using (29) in (23) and by 2-torsion freeness of R, we get

$$h(x)D(x) = 0. \quad (30)$$

and also $D(x)h(x) = 0$. Now take
$h(x)D(y)+2H(x,y)D(x)=h(x)D(y)+ 2([D(x),\alpha(y)]+ [D(y), \alpha(x)])D(x)$.

Replacing y by x and using (9) and (30), we get

$$h(x)D(y) + 2H(x,y)D(x) = 0. \quad (31)$$

$$h(x)D(y) = -2H(x,y)D(x). \quad (32)$$

Replacing y by xy in (31), we get
$h(x)\alpha(x)D(y)+2h(x)\alpha(y)D(x)+2\alpha(x)H(y,x)D(x)+2D(x)[\alpha(y), \alpha(x)]D(x) = 0$.
In view of (32) and (30) the last equation becomes
$[h(x),\alpha(x)]D(y)+2h(x)\alpha(y)D(x)+2D(x)[\alpha(y), \alpha(x)]D(x) = 0$
Using $[h(x),\alpha(x)]= 0$ and 2-torsion freeness of R the last equation becomes

$$h(x)\alpha(y)D(x)+D(x)[\alpha(y),\alpha(x)]D(x)=0. \quad (33)$$

Replacing y by yx in (33), we have

$$h(x)\alpha(y)\alpha(x)D(x)+D(x)[\alpha(y),\alpha(x)]\alpha(x)D(x)=0. \quad (34)$$

Post multiplying (33) by $\alpha(x)$, we get

$$h(x)\alpha(y)D(x)\alpha(x)+D(x)[\alpha(y),\alpha(x)]D(x)\alpha(x)=0. \quad (35)$$

Subtracting (34) from (35) and using (9), we get

$$h(x)\alpha(y)h(x)+D(x)[\alpha(y),\alpha(x)]h(x)= 0. \quad (36)$$

Replacing y by $\alpha^{-1}(t)$ in (36), we get

$$h(x)th(x) + D(x)[t,\alpha(x)]h(x) = 0 \quad (37)$$

Replacing y by $\alpha^{-1}(t)$ in (15) and then replacing z by $\alpha^{-1}(2D(x))$ and using (9), we get
$3h(x)th(x)+4D(x)[t,\alpha(x)]h(x)+2D(x)t[h(x),\alpha(x)]=0$.
As $[h(x),\alpha(x)]=0$, therefore the last equation becomes

$$3h(x)th(x) +4D(x)[t,\alpha(x)]h(x)=0 \quad (38)$$

In view of (37), we get $h(x)th(x) = 0$ and semiprimeness of R implies that $h(x)= 0$ for all $x \in R$. That is $[D(x),\alpha(x)]= 0$ for all $x \in R$. Since $\alpha:R \to R$ is an automorphism, then we obtain
$[D(x), x] = 0$ for all $x \in U$. This lead to $[D(x),x] \in Z(R)$, for all $x \in U$. Thus, according to Lemma 2.1, we get $D(R) \subseteq Z(R)$ and the ideal generated by $D(R)$ is in the $Z(R)$. This completes the proof.

*Corollary 3.7.*

Let R be a 2-torsion free non-commutative semiprime ring and $\alpha:R \to R$ is an automorphism. Suppose that there exists an $(\alpha,\alpha)$-inner derivation $D:R \to R$ defined by $D(x)= \alpha([a,x])$ for all $x \in U$ such that $[[D(x),\alpha(x)],\alpha(x)]= 0$ then D is skew-inner derivation of R.

*Proof:* From (30) in Theorem 3.6, with using R is non-commutative semiprime ring, we obtain $D(x)=0$ for all $x \in U$.
Left-multiplying by x, we obtain
$xD(x)=0$ for all $x \in U$. Again right-multiplying by x, we get $D(x)x=0$ for all $x \in U$. From these equation, we obtain
$D(x)x+xD(x)=0$ for all $x \in U$. Thus, we completed our proof.

By using the similar techniques used in Theorem 2.1 one can prove the following theorem:

*Theorem 3.8.*

Let R be a 2-torsion free and 3-torsion free semiprime ring, U a nonzero left ideal of R such that $r_R(U)=0$, and $\alpha: R \to R$ is an automorphism. If there exists an $(\alpha,\alpha)$-derivation D: $R \to R$ such that $[D(x),\alpha(x)],\alpha(x)] \in Z(R)$, then $D(R) \subseteq Z(R)$ and the ideal generated by D(R) is in the Z(R).

## Acknowledgments

The author thanks the referees for helpful comments and suggestions.